\documentclass[12]{article}
\usepackage{latexsym}
\usepackage{amsfonts}
\usepackage{amsmath}
\usepackage{amssymb}
\usepackage{amsthm}
\usepackage[latin1]{inputenc}

\makeatletter \@addtoreset{equation}{section}

\makeatletter

\begin{document}

\vspace{4cm}
\begin{center} \LARGE{\textbf{Solvable Lie algebras with triangular nilradicals}}
\end{center}
\vspace{1cm}
\large
\begin{center}{S. Tremblay$^{*}$ and P. Winternitz$^{\dagger}$}
\end{center}
\normalsize
\bigskip

\begin{center}
{$*$ Centre de recherches math\'ematiques et D\'epartement de physique, Universit\'e de Montr\'eal, C.P. 6128, succ. Centre-ville, Montr\'eal, QC H3C 3J7, Canada }
\end{center}

\medskip

\begin{center} {$^{\dagger}$ Centre de recherches math\'ematiques et D\'epartement de math\'ematiques et statistiques, Universit\'e de Montr\'eal, C.P. 6128, succ. Centre-ville, Montr\'eal, QC H3C 3J7, Canada}
\end{center}

\medskip

\begin{abstract}
All finite-dimensional indecomposable solvable Lie algebras $L(n,f)$, having the triangular algebra $T(n)$ as their nilradical, are constructed. The number of nonnilpotent elements $f$ in $L(n,f)$ satisfies $1\leq f\leq n-1$ and the dimension of the Lie algebra is $\dim\ L(n,f)=f+\frac{1}{2}n(n-1)$.
\end{abstract}

\normalsize

\section{Introduction}
The purpose of this paper is to construct all indecomposable solvable Lie algebras that
have the ``triangular algebras'' $T(n)$ of dimension $\frac{1}{2}n(n-1)$ ($3\leq n<\infty$) as their
nilradicals. By triangular algebra $T(n)$, we mean the nilpotent Lie algebra isomorphic to the Lie algebra of strictly upper triangular
$n\times n$ matrices. The motivation for such a study is multifold. From a mathematical
point of view, this investigation is part of the classification of all finite dimensional Lie
algebras. The Levi theorem\cite{1,2} tells us that every finite dimensional Lie
algebra $L$ is a semi-direct sum of a semisimple Lie algebra $S$ and a solvable ideal
(the radical $R$):
\begin{equation}
\begin{array}{llll}
L=S\rhd R, & [S,S]=S,& [S,R]\subseteq R,& [R,R]\subset R.
\end{array}
\end{equation}
Semi-simple algebras over fields of complex or real numbers have been classified by
Cartan\cite{3}. However, the classification of solvable Lie algebras is only complete for low
dimensions ($\dim\ L\leq 6$)\cite{4,5,6,7}. From Maltsev\cite{8} we know some important results on
the structure of Lie algebras, but not on solvable Lie algebras with a given nilradical. More
recent articles provided a classification of all Lie algebras with Heisenberg
or Abelian nilradicals\cite{9,10}.

As far as physical applications are concerned, we note that solvable Lie algebras often occur as
Lie algebras of symmetry groups of differential equations\cite{11}. Group invariant solutions can be
obtained by symmetry reduction, using the subalgebras of the symmetry algebra\cite{12}.
In this procedure an important step is to identify the symmetry algebra and its subalgebras
as abstract Lie algebras. A detailed identification presupposes the existence of a
classification of Lie algebras into isomorphism classes.

In Section~2, we formulate the problem and the general strategy that we will adopt. In Section~3,
 we illustrate the procedure for the particular case $n=4$. Guided by this last section, in Section~4 we present
the general classification for arbitrary $n$.

\section{Formulation of the problem}
\subsection{General concepts}
Let us first recall some definitions and known results on solvable Lie algebras.
A Lie algebra $L$ is solvable, if its derived series $DS$ terminates, i.e.
\begin{equation}
DS=\{L_{0}\equiv L,L_{1}=[L,L],\ldots,L_{k}=[L_{k-1},L_{k-1}]=0\}
\end{equation}
for some $k\geq 0$.

A Lie algebra $L$ is nilpotent, if its central series $CS$ also terminates, i.e.
\begin{equation}
CS=\{L_{(0)}\equiv L,L_{(1)}=[L,L_{(0)}],\ldots,L_{(k)}=[L,L_{(k-1)}]=0\}
\end{equation}
for some $k\geq 0$.

The nilradical $NR(L)$ of a solvable Lie algebra $L$ is the maximal nilpotent ideal of $L$.
This nilradical $NR(L)$ is unique and its dimension satisfies\cite{5}
\begin{equation}
\dim\ NR(L)\geq \frac{1}{2}\ \dim\ L.
\end{equation}

Any solvable Lie algebra $L$ can be written as the algebraic sum of the nilradical $NR(L)$ and
a complementary linear space $F$
\begin{equation}
L=F \dot{+} NR(L).
\end{equation}

A Lie algebra $L$ is decomposable if it can, by a change of basis, be transformed into a
direct sum of two (or more) subalgebras
\begin{eqnarray}
L=L_{1} \oplus L_{2},& [L_{1},L_{2}]=0.
\end{eqnarray}

An element $N$ of a Lie algebra $L$ is nilpotent in $L$ if
\begin{eqnarray}
[\ldots [[X,N],N]\ldots N]=0,& \forall X \in L.
\end{eqnarray}

A set of elements $\{X^{\alpha}\}$ of $L$ is linearly nilindependent if no nontrivial linear
combination of them is nilpotent.

A set of matrices $\{A^{\alpha}\}_{\alpha=1\ldots n}$ is linearly nilindependent if no nontrivial linear
combination of them is a nilpotent matrix, i.e.
\begin{eqnarray}
\left(
\sum_{i=1}^{n} c_{i}A^{i}
\right)^{k}=0,
\end{eqnarray}
for some $k\in \mathbb{Z^{+}}$, implies $c_{i}=0\ \forall i$.

\subsection{Basic structure of the Lie algebra and general \mbox{strategy}}
Let us consider the finite triangular algebra $T(n)$ with $n\geq 3$ over the field of complex, or real numbers
($\mathbb{K}=\mathbb{C}\ or \ \mathbb{R}$). A basis for this algebra is
\begin{equation}
\begin{array}{c}
\{N_{ik}\ \mid\ 1\leq i<k\leq n\} \label{eq:base} \\*[2ex]
(N_{ik})_{ab}=\delta_{ia}\ \delta_{kb},\ \ \ \dim\ T(n)=\frac{1}{2}n(n-1)\equiv r.
\end{array}
\end{equation}
This basis can be represented by the standard basis of the strictly upper triangular $n\times n$ matrices. The
commutation relations are
\begin{equation}
[N_{ik},N_{ab}]=\delta_{ka} N_{ib}-\delta_{bi} N_{ak}.
\label{eq:nn}
\end{equation}

We wish to extend this algebra to an indecomposable solvable Lie algebra $L(n,f)$ of
dimension $\frac{1}{2}n(n-1)+f$ having $T(n)$ as its nilradical. In other words, we wish to add
$f$ further linearly nilindependent elements to $T(n)$. Let us denote them
$\{X^{1},\ldots,X^{f}\}$.

The derived algebra $[L,L]$ of a solvable Lie algebra $L$ is contained in its
nilradical\cite{2}. The commutation relations will have the form
\begin{eqnarray}
[X^{\alpha},N_{ik}] &=& A^{\alpha}_{ik\,,\,pq}\ N_{pq}\label{eq:xn} \\*[2ex]
[X^{\alpha},X^{\beta}] &=& \sigma^{\alpha\beta}_{pq}N_{pq}\label{eq:xx} \\*[2ex]
\nonumber
1\leq\alpha,\beta\leq f\leq r,& A^{\alpha}_{ik\,,\,pq},\
\sigma^{\alpha\beta}_{pq}\in\mathbb{K}.
\end{eqnarray}
(Here and in the rest of the paper, we use the Einstein summation convention over repeated indices, unless stated otherwise.)
The commutation relations (\ref{eq:xn}) can be rewritten as
\begin{equation}
\begin{array}{c}
[X^{\alpha},N]=A^{\alpha}N \\*[2ex]
N\equiv(N_{12}\ N_{23}\ldots N_{(n-1)n}\ N_{13}\ldots N_{(n-2)n}\ldots N_{1n})^{T}\\*[2ex]
A^{\alpha}\in \mathbb{K}^{r\times r},\ \ N\in \mathbb{K}^{r\times 1}.\
\end{array}
\label{eq:xamatrice}
\end{equation}
(The superscript $T$ indicates transposition.) Here, $N$ is a ``column vector'' of basis elements of $NR(L)$ ordered by first taking the elements $N_{i(i+1)}$ next to the diagonal, then $N_{i(i+2)}$
(removed two steps from the diagonal), etc. We shall call the matrices $A^{\alpha}$ ``structure matrices''.

A classification of the Lie algebras $L(n,f)$ thus amounts to a classification of the structure matrices
$A^{\alpha}$ and the constants $\sigma^{\alpha\beta}_{pq}$. The Jacobi identities have to be
respected by the following 3 types of triplets (those with 3 $N$'s are satisfies automatically)
\begin{equation}
\{X^{\alpha},N_{ik},N_{ab}\}\ \ f\geq 1,\ \ \{X^{\alpha},X^{\beta},N_{ik}\}\ \ f\geq 2,\ \ \{X^{\alpha},X^{\beta},X^{\gamma}\}\ \ f\geq 3\ ,
\end{equation}
\[
1\leq i<k\leq n,\ \ 1\leq a<b\leq n,\ \ 1\leq \alpha,\beta,\gamma\leq f.
\]
Which give us respectively the 3 equations
\begin{eqnarray}
\delta_{ka}\,A^{\alpha}_{ib\,,\,pq}N_{pq}-\delta_{bi}\,A^{\alpha}_{ak\,,\,pq}N_{pq}+A^{\alpha}_{ik\,,\,bq}N_{aq} \nonumber \\*[2ex]
-A^{\alpha}_{ik\,,\,pa}N_{pb}+A^{\alpha}_{ab\,,\,pi}N_{pk}-A^{\alpha}_{ab\,,\,kq}N_{iq}=0 \label{eq:rjxnn} \\*[2ex]
[A^{\alpha},A^{\beta}]_{ik\,,\,pq}N_{pq}=\sigma^{\alpha\beta}_{kq}N_{iq}-\sigma^{\alpha\beta}_{pi}N_{pk} \label{eq:rjxxn} \\*[2ex]
\sigma^{\alpha\beta}_{pq}A^{\gamma}_{pq\,,\,ik}+\sigma^{\gamma\alpha}_{pq}A^{\beta}_{pq\,,\,ik}+\sigma^{\beta\gamma}_{pq}A^{\alpha}_{pq\,,\,ik}=0.\label{eq:rjxxx}
\end{eqnarray}

The eq.(\ref{eq:rjxnn}) will give restrictions on the form of the structure matrices
$A^{\alpha}$. We will transform these matrices into a ``canonical'' form by transformations that
leave the commutation relations (\ref{eq:nn}) of the $NR(L)$ invariant, but transform the matrices $A^{\alpha}$
and the constants $\sigma^{\alpha\beta}_{pq}$. These transformations will be
\\
\\
(i) Redefinition of all the nonnilpotent elements:
\begin{equation}
\begin{array}{ccl}
X^{\alpha} & \longrightarrow & X^{\alpha}+\mu^{\alpha}_{pq}N_{pq}\ ,\ \ \ \ \ \ \ \ \mu_{pq}^{\alpha}\in
\mathbb{K} \\*[2ex]
\Rightarrow A^{\alpha}_{ik\,,\,ab}&\longrightarrow &
A^{\alpha}_{ik\,,\,ab}+\delta_{kb}\,\mu^{\alpha}_{ai}-\delta_{ia}\,\mu^{\alpha}_{kb}.
\end{array}
\label{eq:transfa1}
\end{equation}
\\
(ii) Change of basis in $NR(L)$:
\begin{equation}
\begin{array}{ccl}
N & \longrightarrow & GN,\ \ \ \ \ \ G\in GL(r,\mathbb{K}) \\*[2ex]
\Rightarrow A^{\alpha} &\longrightarrow& GA^{\alpha}G^{-1}.
\end{array}
\label{eq:transfa2}
\end{equation}
\\
(iii) Linear combinations of the elements $X^{\alpha}$ and hence of the matrices $A^{\alpha}$.
\\
\\
Note that the element $N_{1n}$ does not contribute in the transformation (\ref{eq:transfa1}) since it commutes with all the elements in
the $NR(L)$. Thus $\mu^{\alpha}_{1n}$ is not used in (\ref{eq:transfa1}). Also, the matrix $G$ has to be suitably restricted in order to preserve
the commutation relations (\ref{eq:nn}) of the $NR(L)$.

From the eq.'s (\ref{eq:rjxxn}),(\ref{eq:rjxxx}), we obtain some relations
between the matrices $A^{\alpha}$ and the constants $\sigma^{\alpha\beta}_{pq}$. By exploiting the fact that $\mu^{\alpha}_{1n}$ is not utilized
in (\ref{eq:transfa1}), for $f\geq 2$ we will make an additional transformation to simplify the constants $\sigma^{\alpha\beta}_{pq}$, i.e.
\begin{equation}
\begin{array}{ccl}
X^{\alpha} &\longrightarrow& X^{\alpha}+\mu^{\alpha}_{1n}N_{1n},\ \ \ \ \ \ \ \ \mu^{\alpha}_{1n}\in\mathbb{K} \\*[2ex]
\Rightarrow
\sigma^{\alpha\beta}_{pq} &\longrightarrow&
\sigma^{\alpha\beta}_{pq}+\mu^{\beta}_{1n}\,A^{\alpha}_{1n\,,\,pq}-\mu^{\alpha}_{1n}\,A^{\beta}_{1n\,,\,pq}.
\end{array}
\label{eq:transfsigma}
\end{equation}
(Where $A^{\alpha}\longrightarrow A^{\alpha}$, i.e. the structure matrices stay the same.) It will therefore be possible to simplify some constants
$\sigma^{\alpha\beta}_{pq}$ associated with the matrices $A^{\alpha},A^{\beta}$.

\section{Illustration of the procedure for low\\ dimensions}
\subsection{The case $\mathbf{n=3}$}

In this case, the Lie algebra $T(3)$ is isomorphic to the Heisenberg algebra $H(1)$. As mentioned
previously, solvable Lie algebras with Heisenberg nilradicals were classified earlier \cite{9}. We will therefore
consider $n>3$ from this point on. The dimension $n=3$ is the only case for which there is an
isomorphism between the triangular and the Heisenberg Lie algebras.

\subsection{The case $\mathbf{n=4}$}

In this particular case, we have
\begin{eqnarray}
A^{\alpha}\in \mathbb{K}^{6\times 6},& N=(N_{12}\ N_{23}\ N_{34}\ N_{13}\ N_{24}\ N_{14})^T.
\end{eqnarray}

Let us first consider relations (\ref{eq:rjxnn}). We can separate them into two classes of equations. The first arises from the triplets
$\{X^{\alpha},N_{ik},N_{kb}\},\ 1\leq i<k=a<b\leq 4$, which give
\begin{equation}
\begin{array}{c}
A^{\alpha}_{ib\,,\,pq}N_{pq}+A^{\alpha}_{ik\,,\,bq}N_{kq}-A^{\alpha}_{ik\,,\,pk}N_{pb}
+A^{\alpha}_{kb\,,\,pi}N_{pk}-A^{\alpha}_{kb\,,\,kq}N_{iq}=0\\*[2ex]
(no\ summation\ over\ k).
\end{array}
\label{eq:rjxnn1}
\end{equation}
The second class comes from the triplets $\{X^{\alpha},N_{ik},N_{ab}\},\ 1\leq i<k\leq 4,\ 1\leq a<b\leq 4,\
k\neq a\ (b\neq i)$ and in this case eq.(\ref{eq:rjxnn}) becomes
\begin{equation}
\begin{array}{l}
A^{\alpha}_{ik\,,\,bq}N_{aq}-A^{\alpha}_{ik\,,\,pa}N_{pb}+A^{\alpha}_{ab\,,\,pi}N_{pk}-A^{\alpha}_{ab\,,\,kq}N_{iq}
=0.
\end{array}
\label{eq:rjxnn2}
\end{equation}

We begin by considering eq.(\ref{eq:rjxnn2}). From each possible triplet associated with this class of equation, we use
the linear independence of the $\{N_{lm}\}$ to determine relations between the elements of $A^{\alpha}$. For example, from the triplet $\{X^{\alpha},N_{12},N_{34}\}$, we obtain
\begin{equation}
A^{\alpha}_{12,13}+A^{\alpha}_{34,24}=0,\ \ A^{\alpha}_{12,23}=A^{\alpha}_{34,23}=0.
\end{equation}
When we apply eq.(\ref{eq:rjxnn2}) to the 11 triplets associated to this equation, we find
\begin{equation}
A^{\alpha}=
\left(
\begin{array}{cccccc}
\ast & 0 & A^{\alpha}_{12,34} & A^{\alpha}_{12,13} & \ast & \ast \\
0 & \ast & 0 & \ast & \ast & \ast \\
A^{\alpha}_{34,12} & 0 & \ast & \ast & -(A^{\alpha}_{12,13}) & \ast \\
0 & 0 & 0 & \ast & (A^{\alpha}_{12,34}) & \ast \\
0 & 0 & 0 & (A^{\alpha}_{34,12}) & \ast & \ast \\
0 & 0 & 0 & 0 & 0 & \ast
\end{array}
\right).
\end{equation}
Where $\ast$ denote arbitrary elements unrelated to others in the matrices $A^{\alpha}$.

In the same manner, we apply eq.(\ref{eq:rjxnn1}) to the 4 triplets associated to this class of equation.
This gives us some further relations between the matrix elements and $A^{\alpha}$ becomes
\begin{equation}
A^{\alpha}=
\left(
\begin{array}{cccccc}
A^{\alpha}_{12,12} & 0 & 0 & A^{\alpha}_{12,13} & \ast & \ast \\
& A^{\alpha}_{23,23} & 0 & A^{\alpha}_{23,13} & A^{\alpha}_{23,24} & \ast \\
&& A^{\alpha}_{34,34} & \ast & -(A^{\alpha}_{12,13}) & \ast \\
&&& A^{\alpha}_{12,12}+A^{\alpha}_{23,23} & 0 & (A^{\alpha}_{23,24}) \\
&&&& A^{\alpha}_{23,23}+A^{\alpha}_{34,34} & (A^{\alpha}_{23,13}) \\
&&&&& A^{\alpha}_{12,12}+A^{\alpha}_{23,23}+A^{\alpha}_{34,34}
\end{array}
\right).
\label{eq:a}
\end{equation}

To simplify the form of the matrix (\ref{eq:a}), we carry out the transformation (\ref{eq:transfa1}) for the $f$
matrices $A^{\alpha}$. Given the liberty of the 5 constants $\mu^{\alpha}_{pq}$ for each $\alpha$ independently (the sixth one, $\mu^{\alpha}_{14}$ does not contribute), we can arrange to have
\begin{equation}
A^{\alpha}_{12,13}=A^{\alpha}_{12,14}=A^{\alpha}_{23,13}=A^{\alpha}_{23,24}=A^{\alpha}_{34,14}=0.
\end{equation}
Therefore, each matrix $A^{\alpha}$ can be transformed into
\begin{equation}
\begin{array}{c}
A^{\alpha}=
\left(
\begin{array}{cccccc}
A^{\alpha}_{12,12} &0&0&0& A^{\alpha}_{12,24} & 0\\
& A^{\alpha}_{23,23} &0&0&0& A^{\alpha}_{23,14} \\
&& A^{\alpha}_{34,34} & A^{\alpha}_{34,13} &0&0 \\
&&& A^{\alpha}_{13,13} &0&0 \\
&&&& A^{\alpha}_{24,24} &0 \\
&&&&& A^{\alpha}_{14,14}
\end{array}
\right) \\
\\
A^{\alpha}_{ik\,,\,ik}=\sum_{p=i}^{k-1}A^{\alpha}_{p(p+1)\,,\,p(p+1)}
\end{array}
\label{eq:A}
\end{equation}
These matrices must be linearly nilindependent otherwise the $NR(L)$ would be larger than $T(4)$. In particular, this implies that we cannot
simultaneously have $A^{\alpha}_{12,12}=A^{\alpha}_{23,23}=A^{\alpha}_{34,34}=0$. Also, since we have 3 parameters on the diagonal, the nilindependence between the $A^{\alpha}$ implies that we have at most 3 non-nilpotent elements, i.e.
\begin{equation}
1\leq f\leq 3.
\label{eq:1f3}
\end{equation}

Let us now look at the cases $f\geq 2$. The structure matrices $A^{\alpha}$ have the ``canonical'' form given by (\ref{eq:A}), therefore the possibly nonzero elements of the commutators $[A^{\alpha},A^{\beta}]$ are
\begin{equation}
\begin{array}{lll}
[A^{\alpha},A^{\beta}]_{12,24},& [A^{\alpha},A^{\beta}]_{23,14},& [A^{\alpha},A^{\beta}]_{34,13}.
\end{array}
\label{eq:rcdiff0}
\end{equation}
By the linear independence of the $\{N_{lm}\}$ and from (\ref{eq:rjxxn}),(\ref{eq:rcdiff0}) and (\ref{eq:xx}) we find that
\begin{eqnarray}
[A^{\alpha},A^{\beta}]&=&0 \label{eq:AA} \\*[2ex]
[X^{\alpha},X^{\beta}]&=&\sigma^{\alpha\beta}N_{14}. \label{eq:rcxx}
\end{eqnarray}

Finally, we consider the case $f=3$. In view of the  ``canonical'' form of the structure matrices
$A^{\alpha}$ and by relation (\ref{eq:rcxx}), eq.(\ref{eq:rjxxx}) becomes
\begin{equation}
\sigma^{12}A^{3}_{14,14}+\sigma^{31}A^{2}_{14,14}+\sigma^{23}A^{1}_{14,14}=0. \label{eq:rcxxx-4}
\end{equation}

Moreover, the transformation (\ref{eq:transfsigma}) will modify the constants $\sigma^{\alpha\beta}$ into
\begin{equation}
\sigma^{\alpha\beta}\longrightarrow \sigma^{\alpha\beta}+\mu^{\beta}_{14}\,A^{\alpha}_{14,14}
-\mu^{\alpha}_{14}\,A^{\beta}_{14,14}.
\label{eq:transfxx}
\end{equation}
Hence, by using (\ref{eq:transfxx}) for $f=2$ and (\ref{eq:transfxx}),(\ref{eq:rcxxx-4}) for $f=3$, we obtain
\begin{equation}
[X^{\alpha},X^{\beta}]=
\left\{
\begin{array}{lll}
\sigma^{\alpha\beta} N_{14} &for\ \ A^{1}_{14,14}=\ldots=A^{f}_{14,14}=0 \\
\\
0 &otherwise.
\end{array}
\right.
\label{eq:x1x2}
\end{equation}

To further simplify the structure matrix, let us perform a transformation of the type (\ref{eq:transfa2}), i.e.
\begin{eqnarray}
N\longrightarrow G_{1}N ,&
G_{1}=
\left(
\begin{array}{cccccc}
1 &0&0&0& g_{1} &0 \\
& 1 &0&0&0& g_{2} \\
&& 1 & g_{3} &0&0 \\
&&& 1 &0&0 \\
&&&& 1 &0 \\
&&&&& 1
\end{array}
\right).
\label{eq:G1-4}
\end{eqnarray}
This transformation leaves the commutation relations (\ref{eq:nn}) of the $NR(L)$ invariant, but transforms the matrices
$A^{\alpha}\longrightarrow G_{1}A^{\alpha}G_{1}^{-1}\ \ \forall\alpha$, i.e.
\begin{eqnarray}
A^{\alpha}_{12,24}&\longrightarrow &A_{12,24}^{\alpha}+g_{1}(A^{\alpha}_{23,23}+A^{\alpha}_{34,34}-A^{\alpha}_{12,12})\nonumber\\*[2ex]
A^{\alpha}_{23,14}&\longrightarrow &A_{23,14}^{\alpha}+g_{2}(A^{\alpha}_{12,12}+A^{\alpha}_{34,34})\label{eq:0lambda}\\*[2ex]
A^{\alpha}_{34,13}&\longrightarrow &A^{\alpha}_{34,13}+g_{3}(A^{\alpha}_{12,12}+A^{\alpha}_{23,23}-A^{\alpha}_{34,34}).\nonumber
\end{eqnarray}
Thus using $g_{1}$ we can eliminate $A_{12,24}^{\nu}\ (1\leq \nu\leq f)$ of the specific matrix $A^{\nu}$, provided that
$A^{\nu}_{23,23}+A^{\nu}_{34,34}\neq A^{\nu}_{12,12}$. The constants $g_{2}$ and $g_{3}$ are used in the same way. Therefore, at most 3 off-diagonal elements can be eliminated by this transformation.

Now, we carry out a second transformation $G_{2}$ of the $NR(L)$, such that the total
transformation will be given by $G=G_{2}G_{1}$. The matrix $G_{2}$ is diagonal and the commutation relations (\ref{eq:nn}) of the $NR(L)$ are left invariant for a transformation of the type
\begin{equation}
N\longrightarrow G_{2}N,\
G_{2}=
\left(
\begin{array}{cccccc}
g_{12} &&&&& \\
& g_{23} &&&& \\
&& g_{34} &&& \\
&&& g_{12}g_{23} && \\
&&&& g_{23}g_{34}& \\
&&&&& g_{12}g_{23}g_{34}
\end{array}
\right)\ \ g_{ik}\in \mathbb{K}\!\setminus\!\!\{0\}.
\label{eq:G2-4}
\end{equation}
The matrices $A^{\alpha}$ are transformed as $A^{\alpha}\longrightarrow
G_{2}A^{\alpha}G_{2}^{-1}\  \forall\alpha$, where
\begin{equation}
A^{\alpha}_{12,24} \longrightarrow
\left(
\frac{g_{12}}{g_{24}}
\right)
A^{\alpha}_{12,24},\ \ \
A^{\alpha}_{23,14} \longrightarrow
\left(
\frac{g_{23}}{g_{14}}
\right)
A^{\alpha}_{23,14},\ \ \
A^{\alpha}_{34,13} \longrightarrow
\left(
\frac{g_{34}}{g_{13}}
\right)
A^{\alpha}_{34,13}.
\label{eq:nlambda}
\end{equation}
Therefore, we can normalize up to three non zero off-diagonal elements. For $\mathbb{K}=\mathbb{C}$ they can be set equal to $+1$, for $\mathbb{K}=\mathbb{R}$ we must in some cases allow the possibility of normalizing to either $+1$, or to $-1$.

\subsubsection{The Lie algebras $L(4,1)$}
The matrix $A$ has the ``canonical'' form given by (\ref{eq:A}). Using the transformation (\ref{eq:0lambda}), we can eliminate all
off-diagonal elements, unless the diagonal elements satisfy specific equations (e.g. $A_{23,23}+A_{34,34}-A_{12,12}=0$). At most, two of there equalities can hold, otherwise the matrix $A$ will be nilpotent. Thus, at most two off-diagonal entries remain. They can be normalized to $+1$, with one exception, namely, if we have $A_{12,24}\neq 0,\ A_{34,13}\neq 0$ for $\mathbb{K}=\mathbb{R}$. Then, we can transform to one of the following:
\[A_{12,24}=+1,\ \ \ A_{34,13}=\pm 1\]
($A_{34,13}=-1$ is not equivalent to $A_{34,13}=+1$).

The final result is that for $\mathbb{K}=\mathbb{C}$, 12 inequivalent types of
such matrices exist, one of them depending on 2 complex parameters, 4 depending on one complex parameter, 7 without
parameters. For $\mathbb{K}=\mathbb{R}$, altogether 13 types exist. Among them, one depends on 2 real
parameters, 4 on one real parameter and 8 without parameters. These matrices are listed in Table~A1 in the appendix. The
set of inequivalent matrices~$A$ represent in fact all the possible Lie algebras $L(4,1)$ of dimension~7. For $\mathbb{K}=\mathbb{C}$ the algebra $R_{1,13}$ is equivalent to $K_{1,12}$.

\subsubsection{The Lie algebras $L(4,2)$}
From the eq.s (\ref{eq:AA}),(\ref{eq:x1x2}) and from the previous results on the matrices $A^{\alpha}$, we can now determine the different types of Lie algebras $L(4,2)$ (of dimension 8). For $\mathbb{K}=\mathbb{C}$ or $\mathbb{R}$, 10 inequivalent types
of such algebras exist, one depending on 2 parameters, 5 on one parameter, 4 without parameters. These Lie algebras are presented
in Table~A2 in the appendix.

\subsubsection{The Lie algebra $L(4,3)$}
We can choose a basis for the set of matrices $\{A^{1},A^{2},A^{3}\}$, by putting $A^{\alpha}_{j(j+1)\,,\,j(j+1)}=\delta_{\alpha j}\ (\alpha,j=1,2,3)$ in the general ``canonical'' form. For the matrix $A^{1}$, we use the transformation (\ref{eq:0lambda}) to annul
$A^{1}_{12,24},A^{1}_{23,14}$ and $A^{1}_{34,13}$. The commutativity of the matrices $A^{\alpha}$ imposes that the matrices $A^{2}$ and $A^{3}$ are also diagonal .

Since $A^{1}_{14,14},A^{2}_{14,14},A^{3}_{14,14}$ are different from zero, we can use eq.(\ref{eq:transfxx}) to put
$[X^{1},X^{2}]=[X^{2},X^{3}]=0$ and $[X^{3},X^{1}]=\sigma^{31}N_{14}$. The relation (\ref{eq:rcxxx-4}) then
imposes $\sigma^{31}=0$ and the commutation relations for the non-nilpotent elements become
\begin{equation}
\begin{array}{ll}
[X^{\alpha},X^{\beta}]=0,& \alpha,\beta=1,2,3.
\end{array}
\end{equation}
Therefore, there exists only one Lie algebra $L(4,3)$ ($\mathbb{K}=\mathbb{C}$ or $\mathbb{R}$) ans its dimension is 9. This algebra are is given in Table~A3 in the appendix.

The results of Section~3.2. can be summed up as a theorem.
\newtheorem{theo1}{Theorem}
\begin{theo1}
Every solvable Lie algebra $L(4,f)$ with a six-dimensional triangular nilradical $T(4)$ can be transformed to a canonical basis $\{X^{\alpha},N_{ik}\},\ \alpha=1,\ldots,f,\ 1\leq i<k\leq 4,\ 1\leq f\leq 3\ $. The commutation relations in this basis are given by eq.(\ref{eq:nn}),(\ref{eq:xn}) and (\ref{eq:xx}). The structure matrices $A^{\alpha}$ all have the form (\ref{eq:A}).

For $f=1$ the matrix $A^{1}\equiv A$ has one of the forms given in Table~A1.

For $f=2$ the matrices $\{A^{1},A^{2}\}$ have one of the forms given in Table~A2. The elements $\{X^{1},X^{2}\}$ commute in all cases except $K_{2}$ of Table~A2, when $\sigma$ is a nonzero arbitrary constant.

For $f=3$ there is precisely one such algebra, given by the matrices $\{A^{1},A^{2},A^{3}\}$ of Table~A3, with all elements $\{X^{1},X^{2},X^{3}\}$ commuting.

Every algebra $L(4,f)$ is isomorphic to precisely one algebra in Table~A$f$, for $f=1,2,3$, respectively.
\end{theo1}

\section{Solvable Lie algebras $L(n,f)$ for $n\geq4$}
\subsection{General results}
\newtheorem{lemme1}{Lemma}
\begin{lemme1}
The structure matrices $A^{\alpha}=\{A_{ik\,,\,ab}^{\alpha}\},\ 1\leq i<k\leq n\  ,\ 1\leq a<b\leq n$ have the following properties.
\begin{enumerate}
\item They are upper triangular.
\item The only off-diagonal matrix elements that do not vanish identically and cannot be annuled by a redefinition of the elements $X^{\alpha}$ are:
\begin{equation}
A_{12\,,\,2n}^{\alpha},\ \ A_{j(j+1)\,,\,1n}^{\alpha}\ (2\leq j\leq n-2),\ \ A_{(n-1)n\,,\, 1(n-1)}^{\alpha}.\label{eq:nonzero}
\end{equation}
\item The diagonal elements $A^{\alpha}_{i(i+1)\,,\,i(i+1)},\ 1\leq i\leq n-1$ are free. The other diagonal elements satisfy
\begin{equation}
A^{\alpha}_{ik\,,\,ik}=\sum_{p=i}^{k-1} A^{\alpha}_{p(p+1)\,,\,p(p+1)},\ \ \ \ k>i+1.\label{eq:diag}
\end{equation}
\end{enumerate}
\label{theo:lemmeA}
\end{lemme1}
\emph{Proof.} We shall use relations (\ref{eq:rjxnn}) that are consequences of the Jacobi relations for $\{X^{\alpha},N_{ik},N_{ab}\}$. Let us prove each statement in the theorem separately.

1. All matrix elements below the diagonal vanish identically, i.e.
\begin{equation}
A^{\alpha}_{ik\,,\,ab}=0,\ \ for\
\left\{
\begin{array}{l}
k-i>b-a \\
k-i=b-a,\ i>a.
\end{array}
\right.
\end{equation}
We prove this statement by induction. In Section 3, we have shown that Lemma~\ref{theo:lemmeA} is valid for $n=4$. Now let us assume it is valid for $n=N-1\geq 4$ and prove that it is then also valid for $n=N$. By the induction assumption, we have
\begin{equation}
\begin{array}{c}
A^{\alpha}_{lm\,,\,pq}=0,\ \ for\
\left\{
\begin{array}{l}
m-l>q-p \\
m-l=q-p,\ l>p.
\end{array}
\right. \\
\\
1\leq l<m\leq N-1,\ \ 1\leq p<q\leq N-1.
\end{array}
\end{equation}
Now consider $n=N$. We are adding new entries in old rows $A^{\alpha}_{ik\,,\,aN}$ , new entries in old columns $A^{\alpha}_{iN\,,\,ab}$ and new rows intersecting new columns $A^{\alpha}_{iN\,,\,aN}$ (here, lower case labels vary from $1$ to $N-1$). We must show that all news entries also vanish.

Let us first take (\ref{eq:rjxnn}) for $k=a=i+1,\ 1\leq i\leq N-2,\ i+2\leq b\leq N$. The coefficient of $N_{pq}$ for $p\geq i+2$ provides the identities
\begin{equation}
A^{\alpha}_{ib\,,\,pq}=0.
\label{eq:ncor1}
\end{equation}
In particular we obtain
\begin{equation}
A^{\alpha}_{ib\,,\,pN}=0,\ \ \ \ b-i>N-p,
\end{equation}
which means that we have no nonzero entries in new columns and old rows. Indeed, the smallest possible value of $Z\equiv k-i+a-N$ for which $A^{\alpha}_{ib\,,\,pN}$ is not forced to be zero by eq.(\ref{eq:ncor1}) is reached for $b=i+1,\ p=N-1$ or for $p=i+1,\ b=N-2$. In both cases, the element $A^{\alpha}_{ib\,,\,pN}$ is above the diagonal.

Now consider eq.(\ref{eq:rjxnn}) for $k=a=N-1,\ b=N,\ 1\leq i\leq N-2$. The coefficients of $N_{pq}$ for $q\leq N-2,\ N_{p(N-1)}$ and $N_{pN}$ yield, in particular
\begin{eqnarray}
A^{\alpha}_{iN\,,\,pq}&=&0,\ \ \ \ q\leq N-2 \label{eq:nroc1}\\*[2ex]
A^{\alpha}_{iN\,,\,p(N-1)}&=&0,\ \ \ \ p\geq i \label{eq:nroc2}\\*[2ex]
A^{\alpha}_{iN\,,\,p(N-1)}+A^{\alpha}_{(N-1)N\,,\,pi}&=&0\label{eq:nroc4} \\*[2ex]
A^{\alpha}_{i(N-1)\,,\,p(N-1)}-A^{\alpha}_{iN\,,\,pN}&=&0,\ \ \ \ \ p\neq i. \label{eq:nroc3}
 \end{eqnarray}
Note that eq.(\ref{eq:nroc2}) is obtained from eq.(\ref{eq:nroc4}). We have $A^{\alpha}_{i(N-1)\,,\,p(N-1)}=0$ for $p>i$ by the induction hypothesis. Hence, also $A^{\alpha}_{iN\,,\,pN}=0$ by eq.(\ref{eq:nroc3}). The remaining elements in new rows below the diagonal are $A^{\alpha}_{iN\,,\,(i-1)(N-1)}$ and $A^{\alpha}_{(N-1)N\,,\,(i-1)i}$ with $2\leq i\leq N-1$. Moreover, these elements are related by eq.(\ref{eq:nroc4}) for $2\leq i\leq N-2$. Let us now use relation (\ref{eq:rjxnn}) for $k=i+1,\ a=N-1,\ b=N,\ 1\leq i\leq N-3$. The coefficient of $N_{iq}$ for $q\leq N-1$ must vanish, hence $A^{\alpha}_{(N-1)N\,,\,(i+1)q}=0$ which can be rewritten as
\begin{equation}
A^{\alpha}_{(N-1)N\,,\,(i-1)i},\ \ \ \ 3\leq i\leq N-1.
\end{equation}
The coefficient of $N_{13}$ for $i=2$ must vanish, hence
\begin{equation}
A^{\alpha}_{(N-1)N\,,\,12}=0.
\end{equation}
Relation (\ref{eq:nroc4}) then implies $A^{\alpha}_{iN\,,\,(i-1)(N-1)}=0\ (2\leq i\leq N-2)$ and this completes the proof of the statement that $A^{\alpha}$ is upper triangular.

2. Let us now consider the matrix elements $A^{\alpha}_{ik\,,\,ab}$ above the diagonal. First of all, we note the relations:
\begin{eqnarray}
A^{\alpha}_{ib\,,\,pb}-A^{\alpha}_{ik\,,\,pk}&=&0,\ \ \ \ p\neq i,\ 1\leq i<k<b\leq n \label{eq:ab1} \\*[2ex]
A^{\alpha}_{ib\,,\,iq}-A^{\alpha}_{kb\,,\,kq}&=&0,\ \ \ \ q\neq b ,\ 1\leq i<k<b\leq n\label{eq:ab2} \\*[2ex]
A^{\alpha}_{ik\,,\,ia}+A^{\alpha}_{ab\,,\,kb}&=&0,\ \ \ \ k\neq a,\ b\neq i,\ 1\leq i<k\leq n,\ 1\leq a<b\leq n.\label{eq:ab3}
\end{eqnarray}
Relations (\ref{eq:ab1}) and (\ref{eq:ab2}) follow from (\ref{eq:rjxnn}) with $k=a$, (\ref{eq:ab3}) from (\ref{eq:rjxnn}) with $k\neq a,\ i\neq b$.

Now let us use the transformation (\ref{eq:transfa1}) to annul certain off-diagonal elements. Specifically, we use the coefficient $\mu^{\alpha}_{pq}$ in the following manner:
\begin{equation}
\begin{array}{lll}
\mu^{\alpha}_{1m}: \ \ A^{\alpha}_{m(m+1)\,,\,1(m+1)}&\longrightarrow&0,\ \ \ \ 2\leq m\leq n-1 \\*[2ex]
\mu^{\alpha}_{lm}: \ \ A^{\alpha}_{(l-1)l\,,\,(l-1)m}&\longrightarrow&0,\ \ \ \ 2\leq l\leq n-1,\ l+1\leq m\leq n.
\end{array}
\label{eq:mu}
\end{equation}
Notice that $\mu^{\alpha}_{1n}$ was not used and remains free for future use. Furthermore, combining (\ref{eq:mu}) with (\ref{eq:ab1}),\ldots,(\ref{eq:ab3}) we obtain many more zeros in the matrix $A^{\alpha}$.

Using relations (\ref{eq:rjxnn}) for $k=i+1\neq a,\ b\neq i,\ 1\leq i\leq n-1,\ 1\leq a<b\leq n$ we find the relations
\begin{equation}
\begin{array}{rll}
A^{\alpha}_{i(i+1)\,,\,bq}-A^{\alpha}_{ib\,,\,(i+1)q}&=&0 \\*[2ex]
A^{\alpha}_{ab\,,\,(i+1)q}&=&0,\ \ \ \ i\neq a\neq i+1,\ q\neq b\neq i \\*[2ex]
A^{\alpha}_{i(i+1)\,,\,pa}&=&0,\ \ \ \ a\neq i+1,\ p\neq i \\*[2ex]
A^{\alpha}_{i(i+1)\,,\,bq}&=&0,\ \ \ \ q\neq i+1,\ b\neq i.
\label{eq:ab4}
\end{array}
\end{equation}
We still need information on the elements $A^{\alpha}_{in\,,\,ab},\ A^{\alpha}_{ik\,,\,an}$. For this we consider eq.(\ref{eq:rjxnn}) for $k=a=i+1<b\leq n\ (1\leq i\leq n-2,\ i+2\leq b\leq n)$. We obtain
\begin{equation}
\begin{array}{rll}
A^{\alpha}_{ib\,,\,pq}&=&0,\ \ \ \ i\neq p\neq i+1,\ i+1\neq q\neq b \\*[2ex]
A^{\alpha}_{ib\,,\,p(i+1)}+A^{\alpha}_{(i+1)b\,,\,pi}&=&0.
\label{eq:ab5}
\end{array}
\end{equation}
Together, relations (\ref{eq:ab1}),\ldots,(\ref{eq:ab5}) give us zeros everywhere exept for the elements (\ref{eq:nonzero}). These elements never enter in the eq.(\ref{eq:rjxnn}) exept for some identically respected trivial triplets of the types $\{X^{\alpha},N_{ik},N_{ik}\}$. Therefore, they are all free and this completes the proof of the second affirmation in Lemma~\ref{theo:lemmeA}.

3. To obtain relations between the diagonal elements, take $1\leq i<k=a<b\leq n$ in eq.(\ref{eq:rjxnn}). The coefficient of $N_{ib}$ is
\begin{equation}
A^{\alpha}_{ib\,,\,ib}-A^{\alpha}_{ik\,,\,ik}-A^{\alpha}_{kb\,,\,kb}=0.
\label{eq:eqdiag}
\end{equation}
Choosing $a=i+1,\ b=i+2\ (1\leq i\leq n-2)$, we obtain
\begin{equation}
A^{\alpha}_{i(i+2)\,,\,i(i+2)}=A^{\alpha}_{i(i+1)\,,\,i(i+1)}+A^{\alpha}_{(i+1)(i+2)\,,\,(i+1)(i+2)}=0.
\end{equation}
Now choosing $a=i+2,\ b=i+3\ (1\leq i\leq n-3)$, we obtain
\begin{equation}
A^{\alpha}_{i(i+3)\,,\,i(i+3)}=A^{\alpha}_{i(i+1)\,,\,i(i+1)}+A^{\alpha}_{(i+1)(i+2)\,,\,(i+1)(i+2)}+A^{\alpha}_{(i+2)(i+3)\,,\,(i+2)(i+3)}=0.
\end{equation}
Proceeding recursively, we deduce relation (\ref{eq:diag}) and this completes the proof of statement~ 3 of Lemma~\ref{theo:lemmeA}$\ \Box$
\newtheorem{lemme2}[lemme1]{Lemma}
\begin{lemme2}
The maximal number of nonnilpotent elements is
\begin{equation}
f_{max}=n-1.
\end{equation}
\label{theo:fmax}
\end{lemme2}
\emph{Proof.}
The proof is straightfoward since we have a maximum of $n-1$ parameters on the diagonal and we impose the nilindependence between the matrices $A^{\alpha}.\ \ \Box$

Up to now we have considered the case $f\geq 1$ and this gave us Lemma~\ref{theo:lemmeA} describing each of the matrices $A^{\alpha}$. Now, we shall consider the cases $f\geq 2$ and $f\geq 3$ and must also satisfy the eq.'s (\ref{eq:rjxxn}) and (\ref{eq:rjxxx}).

Let us first consider $f\geq 2$. From Lemma~\ref{theo:lemmeA}, the possibly nonzero elements of the commutators $[A^{\alpha},A^{\beta}]$ are
\begin{equation}
\begin{array}{lll}
[A^{\alpha},A^{\beta}]_{12,2n},& [A^{\alpha},A^{\beta}]_{j(j+1),1n}\ (j=2,\ldots,n-2),& [A^{\alpha},A^{\beta}]_{(n-1)n,1(n-1)}.
\end{array}
\label{eq:rcdiff0-n}
\end{equation}
Therefore, from (\ref{eq:rjxxn}),(\ref{eq:rcdiff0-n}) and (\ref{eq:xx}) we find that
\begin{eqnarray}
[A^{\alpha},A^{\beta}]&=&0 \label{eq:AA-n} \\*[2ex]
[X^{\alpha},X^{\beta}]&=&\sigma^{\alpha\beta}N_{1n}. \label{eq:rcxx-n}
\end{eqnarray}

Finally, we consider $f\geq 3$. From eq.(\ref{eq:rcxx-n}) and Lemma~\ref{theo:lemmeA}, eq.(\ref{eq:rjxxx}) reduces to
\begin{equation}
\sigma^{\alpha\beta}A^{\gamma}_{1n\,,\,1n}+\sigma^{\gamma\alpha}A^{\beta}_{1n\,,\,1n}+\sigma^{\beta\gamma}A^{\alpha}_{1n\,,\,1n}=0.
\label{eq:rjxxxn}
\end{equation}
\newtheorem{lemme3}[lemme1]{Lemma}
\begin{lemme3}
The commutation relations between the structure matrices and the nonnilpotent elements can be transformed to a canonical form satisfying
\begin{eqnarray}
[A^{\alpha},A^{\beta}]&=&0 \\ *[2ex]
[X^{\alpha},X^{\beta}]&=&
\left\{
\begin{array}{ll}
\sigma^{\alpha\beta}N_{1n} & for\ \ A^{1}_{1n\,,\,1n}=\ldots=A^{f}_{1n\,,\,1n}=0 \\*[2ex]
0 & otherwise.
\end{array}
\right.
\label{eq:rcxxmod-n}
\end{eqnarray}
\label{theo:lemmercxx}
\end{lemme3}
\emph{Proof.}
The commutation relations between the structure matrices have been proven already, so we only consider the proof of eq.(\ref{eq:rcxxmod-n}). Using Lemma~\ref{theo:lemmeA} and transformation (\ref{eq:transfsigma}), we modify the constants $\sigma^{\alpha\beta}$ to
\begin{equation}
\sigma^{\alpha\beta}\longrightarrow \sigma^{\alpha\beta}+\mu^{\beta}_{1n}\,A^{\alpha}_{1n\,,\,1n}
-\mu^{\alpha}_{1n}\,A^{\beta}_{1n\,,\,1n}.
\label{eq:transfxx-n}
\end{equation}
Unless we have $A_{1n\,,\,1n}^{1}=\ldots =A_{1n\,,\,1n}^{f}=0$, this transformation can be used to cancel $(f-1)$ constants $\sigma^{\alpha\beta}$. The remaining constants are forced to be zeros by eq.(\ref{eq:rjxxxn}) and this completes the proof.$\ \ \Box$

\subsection{Changes of basis in the nilradical}
As in the case $n=4$, we want to further simplify the structure matrices. For this, we generalize to $n>4$ the previous transformations
$G_{1}$ and $G_{2}$ which transform the $NR(L)$, but preserve its commutation relations. The transformation $G_{1}$ is given by
\begin{equation}
N\longrightarrow G_{1}N,\ \ \ (G_{1})_{ab\,,\,pq}=\delta_{ab\,,\,pq}+\underbrace{\Delta_{ab\,,\,pq}\ g_{a}}_{no\ sum.\ over\ a},\ \ \ \ \ g_{a}\in \mathbb{K}
\end{equation}
where
\begin{equation}
\begin{array}{lll}
\delta_{ab\,,\,pq}&\equiv&\delta_{ap}\ \delta_{bq} \\*[2ex]
\Delta_{ab\,,\,pq}&\equiv&\delta_{ab\,,\,12}\ \delta_{pq\,,\,2n} + \left(\displaystyle\sum_{j=2}^{n-2} \delta_{ab\,,\,j(j+1)}\right) \delta_{pq\,,\,1n} +
\delta_{ab\,,\,(n-1)n}\ \delta_{pq\,,\,1(n-1)}.
\end{array}
\end{equation}
Note that with this definition the elements of $A^{\alpha}$ satisfy $A^{\alpha}_{ab\,,\,pq}=\left(\delta_{ab\,,\,pq}+
\Delta_{ab\,,\,pq}\right) A^{\alpha}_{ab\,,\,pq}$. The transformation preserves the commutation relations in the $NR(L)$ and the matrices
$A^{\alpha}$ are transformed as $A^{\alpha}\longrightarrow G_{1}A^{\alpha}G_{1}^{-1}\ \ \forall\alpha$. The diagonal elements are invariant and
the off-diagonal ones transform as
\begin{equation}
\begin{array}{rll}
A_{12\,,2n}^{\alpha} &\longrightarrow& A_{12\,,\,2n}^{\alpha}+g_{1}\left(A^{\alpha}_{2n\,,\,2n}-A^{\alpha}_{12\,,\,12}\right) \\*[2ex]
A_{j(j+1)\,,\,1n}^{\alpha} &\longrightarrow& A_{j(j+1)\,,\,1n}^{\alpha}+g_{j}\left(A^{\alpha}_{1n\,,\,1n}-A^{\alpha}_{j(j+1)\,,\,j(j+1)}\right),\ \ \ j=2,\ldots,n-2 \\*[2ex]
A_{(n-1)n\,,\,1(n-1)}^{\alpha} &\longrightarrow& A_{(n-1)n\,,\,1(n-1)}^{\alpha}+g_{n-1}\left(A^{\alpha}_{1(n-1)\,,\,1(n-1)}-A^{\alpha}_{(n-1)n\,,\,(n-1)n}\right)
\end{array}
\label{eq:transfA-n}
\end{equation}
with
\[A^{\alpha}_{ik\,,\,ik}=\sum_{p=i}^{k-1}A^{\alpha}_{p(p+1)\,,\,p(p+1)}.\]
As in the case $n=4$, the constants $g_{m}\ (m=1,\ldots,n-1)$ can be used to eliminate up to $(n-1)$ off-diagonal elements of the matrices $A^{\alpha}$.
\newtheorem{lemme4}[lemme1]{Lemma}
\begin{lemme4}
The canonical form of a structure matrix $A^{\alpha}$ has a nonzero off-diagonal element $A^\alpha_{ik\,,\,ab}$ only if
\begin{equation}
A^{\beta}_{ab\,,\,ab}=A^{\beta}_{ik\,,\,ik},\ \ \ \ \beta=1,\ldots,f.
\end{equation}
This is true simultaneously for all $\beta$.
\label{theo:reldiag}
\end{lemme4}
\emph{Proof.} The off-diagonal element $A^{\alpha}_{ik\,,\,ab}$ of a given matrix $A^{\alpha}$ can be transformed to zero by transformation (\ref{eq:transfA-n}), unless we have $A^{\alpha}_{ab\,,\,ab}=A^{\alpha}_{ik\,,\,ik}$. Now let us consider a second matrix $A^{\beta}$. The relation $[A^{\alpha},A^{\beta}]=0$ among the structure matrices implies
\begin{equation}
A^{\alpha}_{ik\,,\,ab}\,(A^{\beta}_{ab\,,\,ab}-A^{\beta}_{ik\,,\,ik})=
A^{\beta}_{ik\,,\,ab}\,(A^{\alpha}_{ab\,,\,ab}-A^{\alpha}_{ik\,,\,ik}),
\end{equation}
so if $A^{\alpha}_{ik\,,\,ab}\neq 0$, we must have $A^{\beta}_{ab\,,\,ab}=A^{\beta}_{ik\,,\,ik}\ \ \forall \beta$.$\ \ \Box$

Now, consider the second transformation $G_{2}$ given by
\begin{equation}
N\longrightarrow G_{2}N,\ \ \ (G_{2})_{ab\,,\,pq}=\delta_{ab\,,\,pq}\ g_{ab},\ \ \ \ \ g_{ab}\in \mathbb{K}\!\setminus\!\!\{0\}.
\end{equation}
This transformation preserves the commutation relations in $NR(L)$ if
\begin{equation}
g_{ab}=\prod_{p=a}^{b-1}g_{p(p+1)}.
\label{eq:gab}
\end{equation}
The matrices $A^{\alpha}$ are transformed as $A^{\alpha}\longrightarrow G_{2}A^{\alpha}G_{2}^{-1}\ \ \forall \alpha$. The off-diagonal elements are transformed as
\begin{equation}
\begin{array}{rll}
A_{12\,,\,2n}^{\alpha} &\longrightarrow&
\left(
\frac{g_{12}}{g_{2n}}
\right)
A_{12\,,\,2n}^{\alpha} \label{eq:epsilon1} \\*[2ex]
A_{j(j+1)\,,\,1n}^{\alpha} &\longrightarrow&
\left(
\frac{g_{j(j+1)}}{g_{1n}}
\right)
A_{j(j+1)\,,\,1n}^{\alpha},\ \ \ \ \ \ j=2,\ldots ,n-2 \\*[2ex]
A_{(n-1)n\,,\,1(n-1)}^{\alpha} &\longrightarrow&
\left(
\frac {g_{(n-1)n}}{g_{1(n-1)}}
\right)
A_{(n-1)n\,,\,1(n-1)}^{\alpha}.
\end{array}
\end{equation}
This transformation is used to normalize the nonzero off-diagonal elements to $+1$ for $\mathbb{K}=\mathbb{C}$ and to
$+1$, or possibly $-1$ for $\mathbb{K}=\mathbb{R}$. Up to $n-1$ elements can be normalized since we have $n-1$ independent entries in $G_{2}$ (see eq.(\ref{eq:gab})).

Note that before the normalization of the off-diagonal elements, we can normalize to $+1$ the first nonzero entry on diagonal of each matrix $A^{\alpha}$ (we choose the first nonzero one).

\subsection{The Lie algebras $L(n,1)$}
Let us consider one extreme case, namely $f=1$, for which we obtain the following lemma.
\newtheorem{lemme5}[lemme1]{Lemma}
\begin{lemme5}
The structure matrix $A=\{A_{ik\,,\,ab}\}\ 1\leq i<k\leq n,\ 1\leq a<b\leq n$ of the Lie algebra $L(n,1)$ has the following properties.
\begin{enumerate}
\item The maximum number of off-diagonal elements is $n-2$.
\item The off-diagonal elements can all be normalized to $+1$ for $\mathbb{K}=\mathbb{C}$ and to $+1$, or $-1$ for $\mathbb{K}=\mathbb{R}$.
\end{enumerate}
\label{theo:ln1}
\end{lemme5}
\emph{Proof.} Let us prove each statement in the theorem separately.

1. First, suppose that we have $n-1$ non-zero off-diagonal elements. The off-diagonal elements remain different from zero if all the terms in the brackets of (\ref{eq:transfA-n}) are zeros. This gives us a system of linear equations for the diagonal elements such that
\begin{equation}
A_{ik\,,\,ik}=0,\ \ \ \ 1\leq i<k\leq n.
\end{equation}
Hence, in this case the matrix $A$ is nilpotent. If we have less or equal to $n-2$ off-diagonal elements, then we obtain at least one free element on the diagonal. In this case the matrix $A$ can (and must) be chosen to be nonnilpotent.

2.Using transformation (\ref{eq:normal-n}) we normalize all $m\leq n-2$ nonzero off-diagonal elements to $+1$. This imposes a system of $n-2$ algebraic constraints on the $n-1$ coefficients $g_{i(i+1)},\ i=1,\ldots, n-1$. These equations always have a solution, but in some cases the solution may be complex. For $\mathbb{K}=\mathbb{C}$ this is consistent. For $\mathbb{K}=\mathbb{R}$ we must modify the initial normalized systems and include the possibility of normalizing to $+1$, or $-1$ the off-diagonal elements. Thus an equivalence class over $\mathbb{C}$ may be split into several over $\mathbb{R}$ (as usual, when restricting from the algebraicly closed field $\mathbb{C}$ to the nonclosed one $\mathbb{R}$).$\ \ \Box$

\subsection{The Lie algebra $L(n,n-1)$}
Let us now consider the other extreme case, namely $f=n-1$.
\newtheorem{lemme6}[lemme1]{Lemma}
\begin{lemme6}
The Lie algebra $L(n,n-1)$ has the following properties.
\begin{enumerate}
\item Only one such Lie algebra exists. The structure matrices $A^{\alpha}$ are all diagonal and can be chosen to satisfy
\begin{equation}
A^{\alpha}_{ik\,,\,ab}=\delta_{ik\,,\,ab}\ \sum_{p=i}^{k-1}\delta_{\alpha p},\ \ \ \ 1\leq i<k\leq n,\ 1\leq \alpha\leq n-1.
\label{eq:An-1}
\end{equation}
\item The nonnilpotent elements always commute, i.e.
\begin{equation}
[X^{\alpha},X^{\beta}]=0,\ \ \ \ 1\leq \alpha,\beta\leq n-1.
\label{eq:xx-lnn-1}
\end{equation}
\end{enumerate}
\label{theo:lnn-1}
\end{lemme6}
\emph{Proof.} Let us prove each statement in the theorem separately.

1. From the Lemma~\ref{theo:lemmeA}, we see that by appropriate linear combinations of the elements $X^{\alpha}$, we can choose $A^{\alpha}$ to satisfy
\begin{equation}
A^{\alpha}_{i(i+1)\,,\,i(i+1)}=\delta_{\alpha i},\ \ \ \ \alpha,i=1,\ldots,n-1.
\end{equation}
By using point 3 of Lemma~\ref{theo:lemmeA}, we obtain (\ref{eq:An-1}). To complete the proof we have to show that the off-diagonal elements are zero for all the structure matrices $A^{\alpha}$. First, we use transformation (\ref{eq:transfA-n}) to annul the off-diagonal elements of $A^{1}$. Commutativity among the structure matrices then implies that the off-diagonal elements of $\{A^{2},\ldots,A^{n-1}\}$ also vanish.

2. The structure matrices given by (\ref{eq:An-1}) satisfy $A^{\alpha}_{1n\,,\,1n}=1\ \forall\alpha$. Hence, from eq.(\ref{eq:rcxxmod-n}) we obtain eq.(\ref{eq:xx-lnn-1}).$\ \ \Box$

\subsection{The main results}
The results of Section~4 constitute the principal result of this article and can be summed up as a theorem.

\newtheorem{theo2}[theo1]{Theorem}
\begin{theo2}
Every solvable Lie algebra $L(n,f)$ with a triangular nilradical $T(n)$ has the dimension $d=f+\frac{1}{2}n(n-1)$ with $1\leq f\leq n-1$. It can be transformed to a canonical basis $\{X^{\alpha},N_{ik}\},\ \alpha=1,\ldots,f,\ 1\leq i<k\leq n\ $with commutation relations
\[
\begin{array}{ccc}
[N_{ik},N_{ab}]=\delta_{ka}N_{ib}-\delta_{bi}N_{ak},&
[X^{\alpha},N_{ik}]=A^{\alpha}_{ik\,,\,pq}N_{pq},&
[X^{\alpha},X^{\beta}]=\sigma^{\alpha\beta}N_{1n}.
\end{array}
\]
The canonical forms of the structure matrices $A^{\alpha}$ and the constants $\sigma^{\alpha\beta}$ satisfy the following conditions:
\begin{enumerate}
\item The matrices $A^{\alpha}$ are linearly nilindependent and have the form specified in Lemma~\ref{theo:lemmeA}. For $f\geq 2$ they all commute, i.e.
$ [A^{\alpha},A^{\beta}]=0.$
\item All constants $\sigma^{\alpha\beta}$ vanish unless we have $A^{\gamma}_{1n\,,\,1n}=0\ $ for $\gamma=1,\ldots,f.$
\item The remaining off-diagonal elements $A^{\alpha}_{ik\,,\,ab}$ also vanish, unless the diagonal elements satisfy $A^{\beta}_{ik\,,\,ik}=A^{\beta}_{ab\,,\,ab}\ $ for $\beta=1,\ldots,f.$
\item When $f$ reaches its maximal value $f=n-1$, then all matrices $A^{\alpha}$ are diagonal as in eq.(\ref{eq:An-1}) and all elements $X^{\alpha}$ commute.
\item For $f=1$ the matrix $A^{1}$ has at most $n-2$ off-diagonal elements that can be normalized as in Lemma~\ref{theo:ln1}.
\end{enumerate}
\label{theo:results}
\end{theo2}
\section{Conclusions}
In this article we have provided a description and classification of solvable Lie algebras with triangular nilradicals $T(n)$. They are nilpotent Lie algebras that are in some sense, the furthest removed from Abelian algebras. Indeed, the dimensions of the Lie algebras in their central series are:
\[\dim\ CS:\ \left(\frac{n(n-1)}{2},\frac{(n-1)(n-2)}{2},\frac{(n-2)(n-3)}{2},\ldots,3,1,0\right).\]

This complements earlier work \cite{9,10} on the classification of solvable Lie algebras with Heisenberg nilradicals (with $\dim\ CS:\ \left(2n+1,1,0\right)$) and Abelian nilradicals (with $\dim\ CS:\ \left(n,0\right)$).

The main results obtained in this paper are summed up in Theorem~\ref{theo:results} of Section~4.

Applications of these algebras are postponed to a forthcoming article. They will concern Lie theory and differential equations. The algebras $L(n,f)$ can appear as symmetry algebras of nonlinear differential equations \cite{11,12}. They will also be used to construct certain nonlinear ordinary differential equations with superposition formulas \cite{13,14,15}.
\section*{Acknowledgements}
The research of P.W. was partially supported by research grants from NSERC of Canada and FCAR du Qu\'ebec.
\section*{Appendix. Lie algebras with a six-dimensional\\ triangular nilradical $T(4)$}
As an illustration of the results obtained above, we give a list of all algebras $L(4,1),\ L(4,2)$ and $L(4,3)$.

We characterize the algebra by the structural matrices $A^{\alpha}$ and by the constants $\sigma^{\alpha\beta}$. For each algebra, we introduce a name $K_{f,i}(a,b)$ or $R_{f,i}(a,b)$. The letter $K$ indicates that the algebra exists both for $\mathbb{K}=\mathbb{C}$ and $\mathbb{K}=\mathbb{R}$; the algebra $R$ is equivalent to some other algebra in the list for $\mathbb{K}=\mathbb{C}$, but inequivalent for $\mathbb{K}=\mathbb{R}$. The first subscript $f$ indicates the number of nonnilpotent elements and the second subscript simply enumerates the algebras. The labels in the brackets indicate parameters in the matrices $A^{\alpha}$. Note that for each matrix we normalize the first nonzero element on the diagonal to $+1$.

\newpage
\scriptsize
\begin{tabular}{lllllllllll}
\textbf{Table A1.} && The Lie algebras $L(4,1)$ \\
\\
\hline
Name && $A$ &&&&&&&& parameters \\
\hline
$K_{1,1}(a,b)$ &&
$\left(
\begin{array}{cccccc}
1 &&&&& \\
& a &&&& \\
&& b &&& \\
&&& 1+a && \\
&&&& a+b & \\
&&&&& 1+a+b
\end{array}
\right)$ &&&&&&&&
$a,b\in \mathbb{K}$ \\
\hline
$K_{1,2}(a)$ &&
$\left(
\begin{array}{cccccc}
0 &&&&& \\
& 1 &&&& \\
&& a &&& \\
&&& 1 && \\
&&&& 1+a & \\
&&&&& 1+a
\end{array}
\right)$ &&&&&&&&
$a \in \mathbb{K}$ \\
\hline
$K_{1,3}$ &&
$\left(
\begin{array}{cccccc}
0 &&&&& \\
& 0 &&&& \\
&& 1 &&& \\
&&& 0 && \\
&&&& 1 & \\
&&&&& 1
\end{array}
\right)$ &&&&&&&& \\
\hline
$K_{1,4}(a)$ &&
$\left(
\begin{array}{cccccc}
1 &&&& 1 & \\
& a &&&& \\
&& 1-a &&& \\
&&& 1+a && \\
&&&& 1 & \\
&&&&& 2
\end{array}
\right)$ &&&&&&&&
$a \in \mathbb{K}$ \\
\hline
$K_{1,5}$ &&
$\left(
\begin{array}{cccccc}
0 &&&& 1 & \\
& 1 &&&& \\
&& -1 &&& \\
&&& 1 && \\
&&&& 0 & \\
&&&&& 0
\end{array}
\right)$ &&&&&&&& \\
\hline
$K_{1,6}(a)$ &&
$\left(
\begin{array}{cccccc}
1 &&&&& \\
& a &&&& 1 \\
&& -1 &&& \\
&&& 1+a && \\
&&&& -1+a & \\
&&&&& a
\end{array}
\right)$ &&&&&&&&
$a \in \mathbb{K}$ \\
\hline
$K_{1,7}$ &&
$\left(
\begin{array}{cccccc}
0 &&&&& \\
& 1 &&&&1 \\
&& 0 &&& \\
&&& 1 && \\
&&&& 1 & \\
&&&&& 1
\end{array}
\right)$ &&&&&&&& \\
\hline
$K_{1,8}(a)$ &&
$\left(
\begin{array}{cccccc}
1 &&&&& \\
& a &&&& \\
&& 1+a & 1 && \\
&&& 1+a && \\
&&&& 1+2a & \\
&&&&& 2(1+a)
\end{array}
\right)$ &&&&&&&&
$a \in \mathbb{K}$ \\
\hline
$K_{1,9}$ &&
$\left(
\begin{array}{cccccc}
0 &&&&& \\
& 1 &&&& \\
&& 1 & 1 && \\
&&& 1 && \\
&&&& 2 & \\
&&&&& 2
\end{array}
\right)$ &&&&&&&& \\
\hline
$K_{1,10}$ &&
$\left(
\begin{array}{cccccc}
1 &&&&& \\
& -2 &&&& 1 \\
&& -1 & 1 && \\
&&& -1 && \\
&&&& -3 & \\
&&&&& -2
\end{array}
\right)$ &&&&&&&& \\
\hline
\end{tabular}

\begin{tabular}{llllll}
\textbf{Table A1.} && (continued) \\
\\
\hline
Name && $A$ && parameters \\
\hline
$K_{1,11}$ &&
$\left(
\begin{array}{cccccc}
1 &&&& 1 & \\
& 2 &&&& 1 \\
&& -1 &&& \\
&&& 3 && \\
&&&& 1 & \\
&&&&& 2
\end{array}
\right)$ && \\
\hline
$K_{1,12}$ &&
$\left(
\begin{array}{cccccc}
1 &&&& 1 & \\
& 0 &&&& \\
&& 1 & 1 && \\
&&& 1 && \\
&&&& 1 & \\
&&&&& 2
\end{array}
\right)$ && \\
\hline
$R_{1,13}$ &&
$\left(
\begin{array}{cccccc}
1 &&&& 1 & \\
& 0 &&&& \\
&& 1 & -1 && \\
&&& 1 && \\
&&&& 1 & \\
&&&&& 2
\end{array}
\right)$ && \\
\hline
\\
\\
\\
\textbf{Table A2.} && The Lie algebras $L(4,2)$ \\
\\
\hline
Name & $\sigma$ & $A^{1}$ & $A^{2}$ & parameters \\
\hline
$K_{2,1}(a,b)$ & 0 &
$\left(
\begin{array}{cccccc}
1 &&&&& \\
& 0 &&&& \\
&& a &&& \\
&&& 1 && \\
&&&& a & \\
&&&&& 1+a
\end{array}
\right)$ &
$\left(
\begin{array}{cccccc}
0 &&&&& \\
& 1 &&&& \\
&& b &&& \\
&&& 1 && \\
&&&& 1+b \\
&&&&& 1+b
\end{array}
\right)$ &
$a,b \in \mathbb{K}$ \\
\hline
$K_{2,2}$ &
$\sigma$ &
$\left(
\begin{array}{cccccc}
1 &&&&& \\
& 0 &&&& \\
&& -1 &&& \\
&&& 1 && \\
&&&& -1 & \\
&&&&& 0
\end{array}
\right)$ &
$\left(
\begin{array}{cccccc}
0 &&&&& \\
& 1 &&&& \\
&& -1 &&& \\
&&& 1 && \\
&&&& 0 & \\
&&&&& 0
\end{array}
\right)$ &
$\sigma \in \mathbb{K}\!\setminus\!\{0\}$ \\
\hline
$K_{2,3}(a)$ & 0 &
$\left(
\begin{array}{cccccc}
1 &&&&& \\
& a &&&& \\
&& 0 &&& \\
&&& 1+a && \\
&&&& a & \\
&&&&& 1+a
\end{array}
\right)$ &
$\left(
\begin{array}{cccccc}
0 &&&&& \\
& 0 &&&& \\
&& 1 &&& \\
&&& 0 && \\
&&&& 1 & \\
&&&&& 1
\end{array}
\right)$ &
$a\in \mathbb{K}$ \\
\hline
$K_{2,4}$ & 0 &
$\left(
\begin{array}{cccccc}
0 &&&&& \\
& 0 &&&& \\
&& 1 &&& \\
&&& 0 && \\
&&&& 1 & \\
&&&&& 1
\end{array}
\right)$ &
$\left(
\begin{array}{cccccc}
0 &&&&& \\
& 1 &&&& \\
&& 0 &&& \\
&&& 1 && \\
&&&& 1 & \\
&&&&& 1
\end{array}
\right)$ \\
\hline
$K_{2,5}(a)$ & 0 &
$\left(
\begin{array}{cccccc}
1 &&&&& \\
& a &&&& \\
&& 1-a &&& \\
&&& 1+a && \\
&&&& 1 & \\
&&&&& 2
\end{array}
\right)$ &
$\left(
\begin{array}{cccccc}
0 &&&& 1 & \\
& 1 &&&& \\
&& -1 &&& \\
&&& 1 && \\
&&&& 0 & \\
&&&&& 0
\end{array}
\right)$ &
$a \in \mathbb{K}$ \\
\hline
$K_{2,6}$ & 0 &
$\left(
\begin{array}{cccccc}
0 &&&&& \\
& 1 &&&& \\
&& -1 &&& \\
&&& 1 && \\
&&&& 0 & \\
&&&&& 0
\end{array}
\right)$ &
$\left(
\begin{array}{cccccc}
1 &&&& 1 & \\
& 0 &&&& \\
&& 1 &&& \\
&&& 1 && \\
&&&& 1 & \\
&&&&& 2
\end{array}
\right)$ & \\
\hline
\end{tabular}

\begin{tabular}{lllll}
\textbf{Table A2.} && (continued) \\
\\
\hline
Name & $\sigma$ & $A^{1}$ & $A^{2}$ & parameters \\
\hline
$K_{2,7}(a)$ & 0 &
$\left(
\begin{array}{cccccc}
1 &&&&& \\
& a &&&& \\
&& -1 &&& \\
&&& 1+a && \\
&&&& -1+a & \\
&&&&& a
\end{array}
\right)$ &
$\left(
\begin{array}{cccccc}
0 &&&&& \\
& 1 &&&& 1 \\
&& 0 &&& \\
&&& 1 && \\
&&&& 1 & \\
&&&&& 1
\end{array}
\right)$ &
$a \in \mathbb{K}$ \\
\hline
$K_{2,8}$ & 0 &
$\left(
\begin{array}{cccccc}
0 &&&&& \\
& 1 &&&& \\
&& 0 &&& \\
&&& 1 && \\
&&&& 1 & \\
&&&&& 1
\end{array}
\right)$ &
$\left(
\begin{array}{cccccc}
1 &&&&& \\
& 0 &&&& 1 \\
&& -1 &&& \\
&&& 1 && \\
&&&& -1 & \\
&&&&& 0
\end{array}
\right)$ & \\
\hline
$K_{2,9}(a)$ & 0 &
$\left(
\begin{array}{cccccc}
1 &&&&& \\
& a &&&& \\
&& 1+a &&& \\
&&& 1+a && \\
&&&& 1+2a & \\
&&&&& 2(1+a)
\end{array}
\right)$ &
$\left(
\begin{array}{cccccc}
0 &&&& \\
& 1 &&&& \\
&& 1 & 1 && \\
&&& 1 && \\
&&&& 2 &\\
&&&&& 2
\end{array}
\right)$ &
$a\in \mathbb{K}$ \\
\hline
$K_{2,10}$ & 0 &
$\left(
\begin{array}{cccccc}
0 &&&&& \\
& 1 &&&& \\
&& 1 &&& \\
&&& 1 && \\
&&&& 2 & \\
&&&&& 2
\end{array}
\right)$ &
$\left(
\begin{array}{cccccc}
1 &&&&& \\
& 0 &&&& \\
&& 1 & 1 && \\
&&& 1 && \\
&&&& 1 & \\
&&&&& 2
\end{array}
\right)$ & \\
\hline
\\
\\
\\
\end{tabular}
\begin{tabular}{lllllllllll}
\textbf{Table A3.} && The Lie algebra $L(4,3)$ &&&&&& \\
\\
\hline
Name & $\sigma's$ & $A^{1}$ & $A^{2}$ & $A^{3}$ \\
\hline
$K_{3,1}$ & 0 &
$\left(
\begin{array}{cccccc}
1 &&&&& \\
& 0 &&&& \\
&& 0 &&& \\
&&& 1 && \\
&&&& 0 & \\
&&&&& 1
\end{array}
\right)$ &
$\left(
\begin{array}{cccccc}
0 &&&&& \\
& 1 &&&& \\
&& 0 &&& \\
&&& 1 && \\
&&&& 1 & \\
&&&&& 1
\end{array}
\right)$ &
$\left(
\begin{array}{cccccc}
0 &&&&& \\
& 0 &&&& \\
&& 1 &&& \\
&&& 0 && \\
&&&& 1 & \\
&&&&& 1
\end{array}
\right)$  \\
\hline
\end{tabular}

\end{document}